\documentclass[10pt,a4paper]{article}
\usepackage[utf8]{inputenc}
\usepackage[english]{babel}

\usepackage{amsmath}
\usepackage{amsfonts}
\usepackage{amssymb}
\usepackage{amsthm}

\usepackage{graphicx}
\usepackage[margin=2cm]{geometry}
\usepackage{paralist}

\usepackage{csquotes,etoolbox,keyval,ifthen,url,babel}
\usepackage[backend=biber, defernumbers=true, sorting=anyt, style=numeric-comp, maxbibnames=4, maxcitenames=2]{biblatex}
\usepackage{hyperref}

\hypersetup{
    colorlinks=true,
    linkcolor=red,
    filecolor=magenta,      
    urlcolor=cyan,
    breaklinks=true
            }

\usepackage{chngcntr}
\usepackage{caption}
\usepackage{subcaption}
\usepackage{wrapfig}

\usepackage{booktabs}

\newtheorem{theorem}{Teorema}

\newtheorem{example}{Primjer}
\newtheorem{remark}{Primjedba}

\counterwithin{theorem}{section}
\counterwithin{definition}{section}
\counterwithin{example}{section}
\counterwithin{lemma}{section}
\counterwithin{remark}{section}

\DeclareMathOperator{\ord}{Ord}

\graphicspath{{Slike/}}

\allowdisplaybreaks
\usepackage{mathtools}
\usepackage[multiple]{footmisc} 

\title{O numeri\v ckom rje\v savanju Cauchyjevog problema Runge-Kutta metodama na Shishkinovoj mre\v zi}
\author{Vesna Divkovi\'c$^{1}$,
	    Laura Luki\'c$^{1}$, 
	    Elvir Memi\'c$^{1}$ i 
	    Samir Karasulji\'c$^{1,\,2}$} 
\date{}

\begin{document}
	\maketitle
	\footnotetext[1]{Prirodno--matemati\v cki fakultet Univerziteta u Tuzli}
	\footnotetext[2]{corresponding author}
	\begin{abstract}
		U ovom radu razmatrano je numeri\v cko rje\v savanje singularno--perturbacionog Cauchyjevog problema Runge--Kutta metodama na Shishkinovoj mre\v zi. Numeri\v cka rje\v senja posmatranog problema dobijena su kori\v stenjem dvije eksplicitne i jedne implicitne Runge--Kutta metode na najjednostavnijoj slojno--adaptivnoj mre\v zi. Na kraju su dobijeni rezultati upore\dj eni.
	\end{abstract}
	
\section{Uvod}	

Diferencijalne jedna\v cine koriste se za modeliranje raznih problema u prirodnim, in\v zinjerskim pa \v cak i u dru\v stvenim naukama. U velikom broju ovih problema zahtijeva se da rje\v senje diferencijalne jedna\v cine zadovoljava i jedan dodatni uslov, po\v cetni uslov ili po\v cetnu vrijednost. 

U realnim problemima, diferencijalne jedna\v cine koje sre\'cemo u matemati\v ckim modelima, su i suvi\v se te\v ske da bi ih ta\v cno rije\v sili, a nekada je to i nemogu\'ce.  Postoje dva pristupa za prevazila\v zenje prethodne situacije. Prvi pristup je u pojednostavljivanju date diferencijalne jedna\v cine ili matemati\v ckog modela, tako da mo\v zemo izra\v cunati ta\v cno rje\v senje diferencijalne jedna\v cine, zatim dobijeno rje\v senje koristimo kao aproksimaciju rje\v senja realnog tj. originalnog problema. Drugi pristup je da odmah ra\v cunamo aproksimativno rje\v senje ili preciznije numeri\v cko rje\v senje. U najve\'cem broju slu\v cajeva, drugi pristup je bolji. Dakle, koristimo bolji matemati\v cki model, koji je "bli\v zi" realnom problemu, te ra\v cunamo odgovaraju\'ce numeri\v cko rje\v senje, jer se skoro po pravilu u ovakvim matemati\v ckim modelima pojavljuju "komplikovanije" diferencijalne jedna\v cine, \v cija ta\v cna rje\v senje ili ne mo\v zemo izra\v cunati ili je to veoma te\v sko.

U nastavku ovog rada razmatra\'cemo sljede\'ci Cauchyjev problem

\begin{equation}\label{problem1}
   \begin{cases}
      y'=f(x,y),\\
      y(x_0)=y_0.
   \end{cases}
\end{equation}
Postavljaju se pitanja kada Cauchyjev problem ima rje\v senje i ako ono postoji da li je to rje\v senje jedinstveno? Odgovori na ova pitanja dati su u sljede\'coj teoremi.
\begin{theorem}{\rm\cite[ pp 447]{cheney2004}}\label{teorema1}
	Ako su $f$ i $\frac{\partial f}{\partial y}$ neprekidni na pravougaoniku definisanom sa $|x-x_0|<\alpha$  i $|y-y_0|<\beta,$ tada Cauchyjev problem \eqref{problem1} ima jedinstveno neprekidno rje\v senje na nekom intervalu $|x-x_0|<\gamma.$
\end{theorem} 
Vrijednost konstante $\gamma$ je najmanje $\frac{\beta}{M},$ gdje je $M$ gornja granica za $|f(x,y)|$ na pravougaoniku definisanom u upravo navedenoj teoremi. 

Kako je ve\'c pomenuto, vrlo je uska klasa diferencijalnih jedna\v cina koje se mogu ta\v cno rije\v sti, pa se stoga pribjegava ra\v cunanju numeri\v ckog rje\v senja. Iz standardnih kurseva numeri\v cke matematike poznate su metode za ra\v cunanje ovakvog rje\v senja, npr. Eulerova i njene modifikacije, Taylorova,  Runge--Kutta, vi\v sekora\v cne i dr. \ \\ \ \\

Osim Cauchyjevog problema \eqref{problem1} \v cesto sre\'cemo i njegovu modifikaciju kod koje je prvi izvod pomno\v zen nekim pozitivnim malim parametrom $\varepsilon,$ ($0<\varepsilon\leqslant 1$ ).  
\begin{equation}\label{problem2a}
\begin{cases}
\varepsilon y'=\tilde{f}(x,y),\quad 0<x,\quad 0< \varepsilon \leqslant 1,\\
y(x_0,\varepsilon)=y_0,
\end{cases}
\end{equation}
gdje je $f(x,y)\in C^{n,n}([0,a]\times \mathbb{R}),\: n\geqslant 1.$ Podijeliv\v si prethodnu diferencijalnu jedna\v cinu  sa parametrom $\varepsilon,$ dobijamo 
\begin{equation}\label{problem2}
\begin{cases}
  y'=f(x,y),\quad 0<x,\quad 0< \varepsilon \leqslant 1,\\
y(x_0,\varepsilon)=y_0,
\end{cases}
\end{equation}
gdje je $f(x,y)=\frac{\tilde{f}(x,y)}{\varepsilon}.$ Neka je $\tilde{f}$ linearna funkcija po $y,$ tj. neka vrijedi  $\tilde{f}(x,y)=\tilde{p}(x)y+\tilde{q}(x),$ u tom slu\v caju Cauchyjev problem \eqref{problem2a} poprima sljede\'ci oblik 
\begin{equation*}
  \begin{cases}
    \varepsilon y'=\tilde{p}(x)y+\tilde{q}(x),
    \quad 0<x,\quad 0< \varepsilon \leqslant 1,\\
    y(x_0, \varepsilon)=y_0.
  \end{cases}
\end{equation*}
Ponovo podijeliv\v si prethodnu diferencijalnu jedna\v cinu sa $\varepsilon,$ dobijamo
\begin{equation}\label{problem3}
   \begin{cases}
       y'=p(x)y+q(x),\quad 
          0<x,\quad 0< \varepsilon \leqslant 1,\\        
          y(x_0,\varepsilon)=y_0,
   \end{cases}
\end{equation}
gdje su $p(x)=\frac{\tilde{p}(x)}{\varepsilon}$ i $q(x)=\frac{\tilde{q}}{\varepsilon}.$ \\

Prisustvo parametra $\varepsilon$ dovodi do brzih promjena ta\v cnog rje\v senja $y$ Cachyjevog problema \eqref{problem2} odnosno \eqref{problem3}  u nekim  dijelovima domena, o \v cemu \'ce biti vi\v se rije\v ci u sljede\'coj sekciji. Zbog ovih brzih promjena, klasi\v cne metode su neadekvatne za  numeri\v cko rje\v savanje navedenih problema u kojima se pojavljuje perturbacioni parametar $\varepsilon$, stoga je bilo potrebno razviti nove efikasnije metode u kojima se uzima u obzir postojanje pomenutih brzih promjena ta\v cnog rje\v senja. 

Jedna od naj\v sirenijih metoda za rje\v savanje ovakvih problema je metoda slojno--adaptivnih mre\v za. Procjena ta\v cnog rje\v senja i njegovih izvoda je veoma va\v zna komponenta u konstruisanju slojno--adaptivnih mre\v za. U sljede\'coj teoremi date su ove procjena za problema \eqref{problem2}, odnosno za ta\v cno rje\v senje ovog problema.

\begin{theorem}{\rm \cite[pp 66]{liseikin2001layer}}\label{teorema2}
	 Neka je $y(x,\varepsilon)$ rje\v senje problema \eqref{problem2}. Tada za $0\leqslant i\leqslant n$ i  $0\leqslant x\leqslant a,$ vrijede sljede\'ce procjene 
	 \begin{equation}\label{procjena1}
	  \left|y^{(i)}(x,\varepsilon)\right|\leqslant C \left\{ 1+\varepsilon^{-i}\exp(-c(0)x/\varepsilon )   \right\},
	 \end{equation}
	 ako je $f_y(x,y)\geqslant c(x)>0$ i $c(x)\in C[0,a];$
	 \begin{equation}
	    \left| y^{(i)}(x,\varepsilon) \right|\leqslant C \left[ 1+\varepsilon^{-i/(k+1)} \exp(-mx^{k+1}/\varepsilon) + (\varepsilon^{1/(k+1)}+x)^{1-i}\right],
	 \end{equation}
     ako je $f_y(x,y)=x^kg(x,y)$ i $g_y(x,y)\geqslant c(x)>0,$ gdje je $c(x)\in C[0,a],$ $0<m<c\min_{x\in[0,a]}c(x),$ $k\geqslant 1$ je pozitivan cio broj. 
	 
\end{theorem}

\begin{remark}
	Posmatrani Cauchyjev problem \eqref{problem1} je najjednostavniji, naime u ovom se problemu pojavljuje diferencijalna jedna\v cina prvog reda. Umjesto diferencijalne jedna\v cine prvog reda, mo\v ze se pojaviti i diferencijalna jedna\v cina vi\v seg reda kao i sistem diferencijalnih jedna\v cina. Svi ovi slu\v cajevi su detaljno obra\dj eni u literaturi. Cilj ovog rada je da uka\v ze na probleme prilikom numeri\v ckog  rje\v savanja Cauchyjevih problema  \v cija ta\v cna rje\v senja imaju brze promjena, pa je sasvim dovoljno u ovoj nekoj po\v cetnoj fazi posmatrati samo Cauchyjeve problema sa diferencijalnom jedna\v cinom prvog reda.  
\end{remark}

\section{Slojno--adaptivne mre\v ze}

U ovoj sekciji dati su osnovni razlozi i ideje koje su dovele do konstruisanja slojno--adaptivnih mre\v za, kao i konstrukcija najjednostavnije slojno--adaptivne mre\v ze.  Analiziran je jednostavan Cauchyjev problem \v cije je ta\v cno rje\v senje poznato. Numeri\v cko rje\v senje za razli\v cite vrijednosti perturbacionog parametra $\varepsilon$ je izra\v cunato na ekvidistantnoj mre\v zi i istaknuti su nedostaci u ovakvom pristupu. Upravo ovi nedostaci doveli su razvoja slojno--adaptivih mre\v za. \\

\begin{figure}[!h]
	\centering 
		\includegraphics[scale=.5]{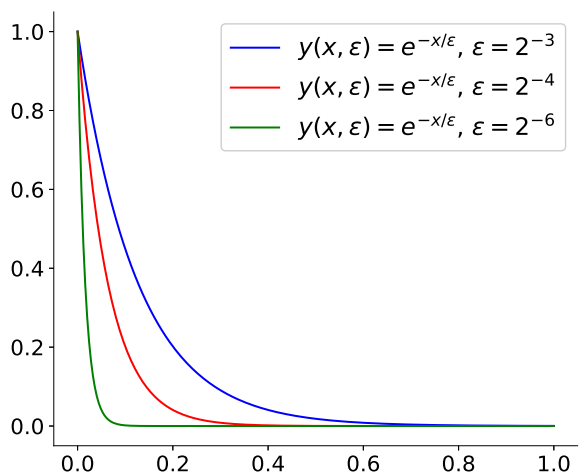}
		\caption{Grafik funkcije $y(x,\varepsilon)=e^{-x/\varepsilon}$ za razli\v cite vrijednosti parametra $\varepsilon$}
		\label{tacno1}
\end{figure}

Posmatrajmo jednostavan testni Cauchyjev problem 
\begin{equation}\label{cauchy1}
\begin{cases}
\varepsilon y'=-y,\: x\in[0,1]\\
y(0,\varepsilon)=1.
\end{cases}
\end{equation}
Ta\v cno rje\v senje ovog problema je $y(x,\varepsilon)=e^{-x/\varepsilon}.$ Grafici funkcije $y(x,\varepsilon)$ za $x\in[0,1]$ dati su na slici \ref{tacno1}, za tri razli\v cite vrijednosti parametra $\varepsilon,\: \varepsilon=2^{-3},\, 2^{-4}$ i $2^{-6}.$ Sa grafika je lako uo\v citi, a to se da zaklju\v citi i na osnovu osobina  funkcije $x\mapsto e^{-x},$ da funkcija ima brze promjene  u nekoj okolini ta\v cke $x=0.$  

Ove promjene su br\v ze \v sto je manji parametar $\varepsilon.$ Ovaj dio domena gdje se de\v savaju ove brze promjene nazivamo sloj i ovo je sloj eksponencijalnog tipa.

Kao \v sto je poznato, u op\v stem slu\v caju Cauchyjev problem ne mo\v zemo ta\v cno rije\v siti, pa se pribjegava numeri\v ckim metodama. U najve\'cem broju tih metoda podijelimo domen, na kojem je potrebno da odredimo rje\v senje, sa odre\dj enim brojem \v cvorova i primjenimo neku od metoda kojom \'cemo izra\v cunati numeri\v cko rje\v senje. Isto tako u najve\'cem broju slu\v caju koristimo ravnomjerno ili uniformnu raspodjelu \v cvorova, tj. rastojanje izme\dj u dva susjedna \v cvora je ekvidistantno. Skup \v cvorova kojim dijelimo neki segment nazivamo mre\v za. Ako npr. na segmentu $[0,1]$  \v zelimo izra\v cunati numeri\v cko rje\v senje nekog Cauchyjevog problema, podijeli\'cemo ga \v cvorovima za koje vrijedi 
\[0=x_0<x_1<\ldots<x_{i-1}<x_i<x_{i+1}<\ldots<x_n=1.\] 
Skup ovih \v cvorova $x_i,\: i=0,1,\ldots,N$ je mre\v za. 
\begin{remark}
U savremenoj literaturi \v ce\v s\'ce se koristi izraz ta\v cka mre\v ze umjesto \v cvor,  pa \'cemo u nastavku ovog rada koristi izraz ta\v cka/ta\v cke mre\v ze.
\end{remark}
Ako su ta\v cke mre\v ze uniformno raspore\dj ene onda je rije\v c je o uniformnoj ili ekvidistantnoj mre\v zi, u suprotnom o neekvidistantnoj ili neuniformnoj mre\v zi. Korak mre\v ze ili parametar mre\v ze je rastojanje izme\dj u dvije  susjedne ta\v cke mre\v ze. Kod uniformnih mre\v za korak mre\v ze je konstantan, obi\v cno ga ozna\v cavamo $h$ i vrijedi $h=\frac{1}{N}$ ako je segment $[0,1]$ na kojem se ra\v cuna numeri\v cko rje\v senje, u slu\v caju segmenta $[a,b]$ vrijedi  $h=\frac{b-a}{N}.$ Kod neuniformnih mre\v za korak mre\v ze ra\v cunamo $h_i=x_{i+1}-x_i,\: i=0,1,\ldots,N-1.$

Uniformne mre\v ze su najjednostavnije pa je samim tim i analiza metoda, koje koriste ovakve mre\v ze, jednostavnija od analize metoda koje koriste neuniformne mre\a v ze. Me\dj utim, u mnogim slu\v cajevima  uniformne mre\v ze nisu najbolje rje\v senje za primjenu. 

\begin{figure}\centering
	\includegraphics[scale=.35]{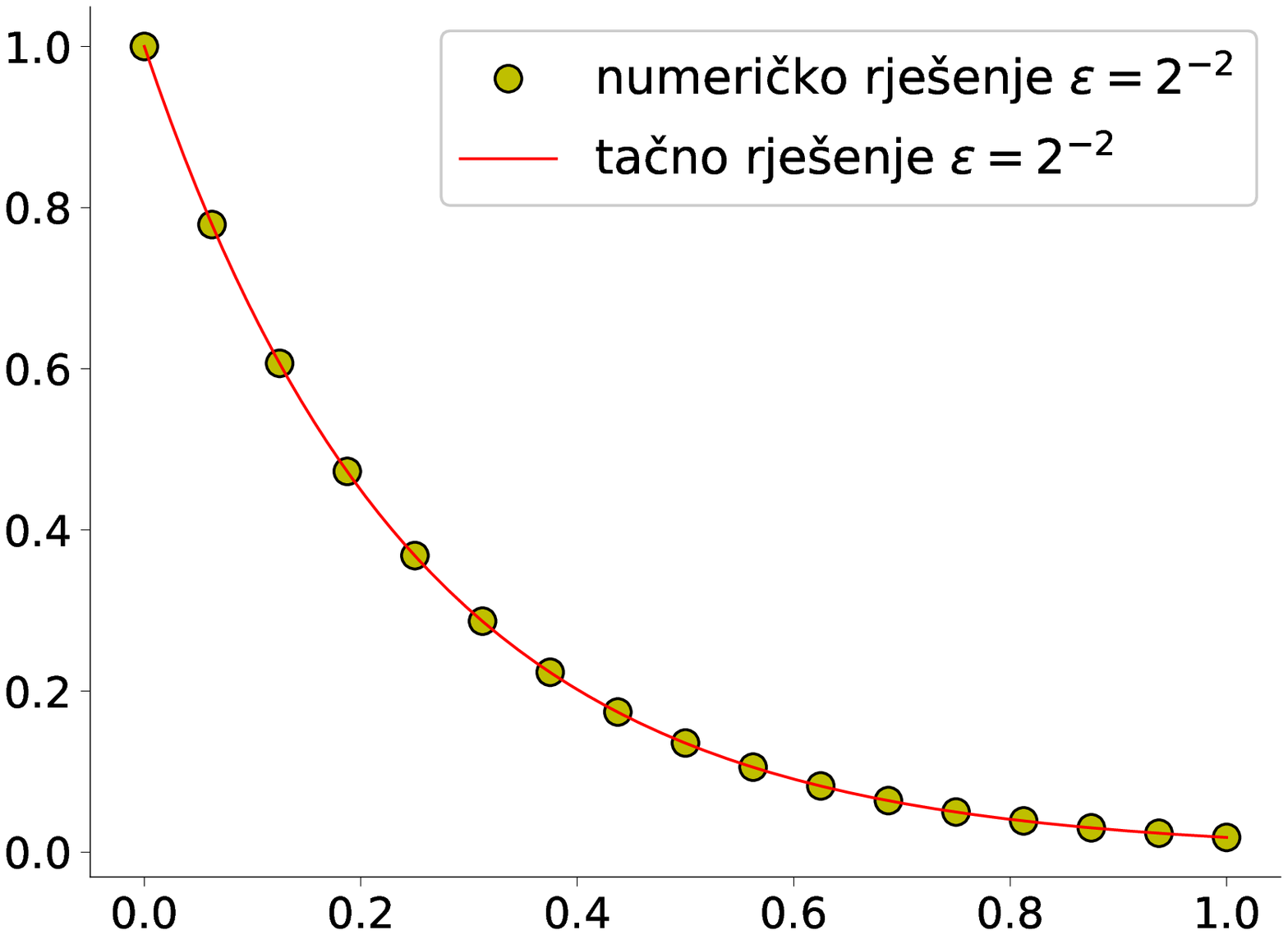}
	\includegraphics[scale=.35]{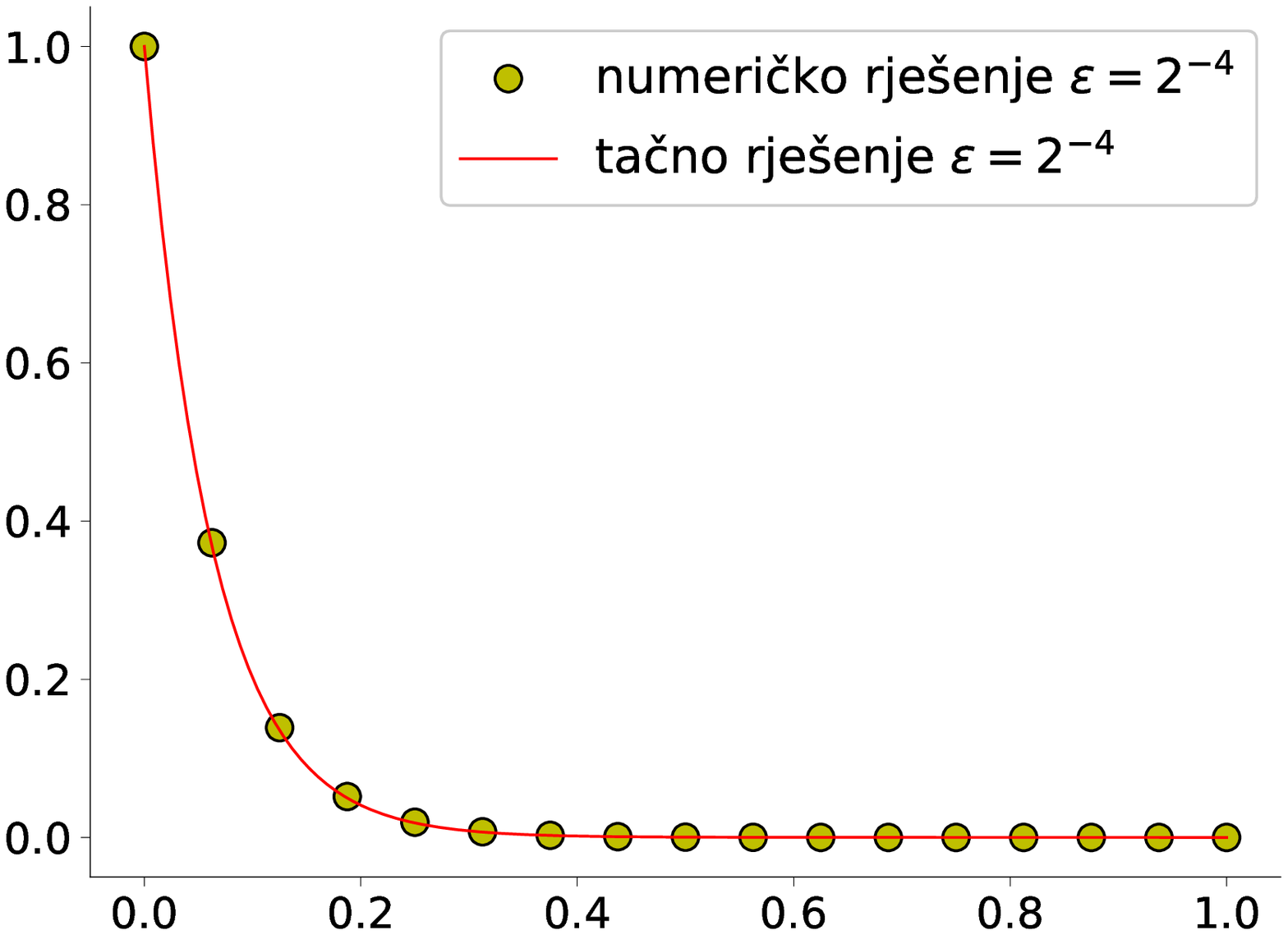}
	\includegraphics[scale=.35]{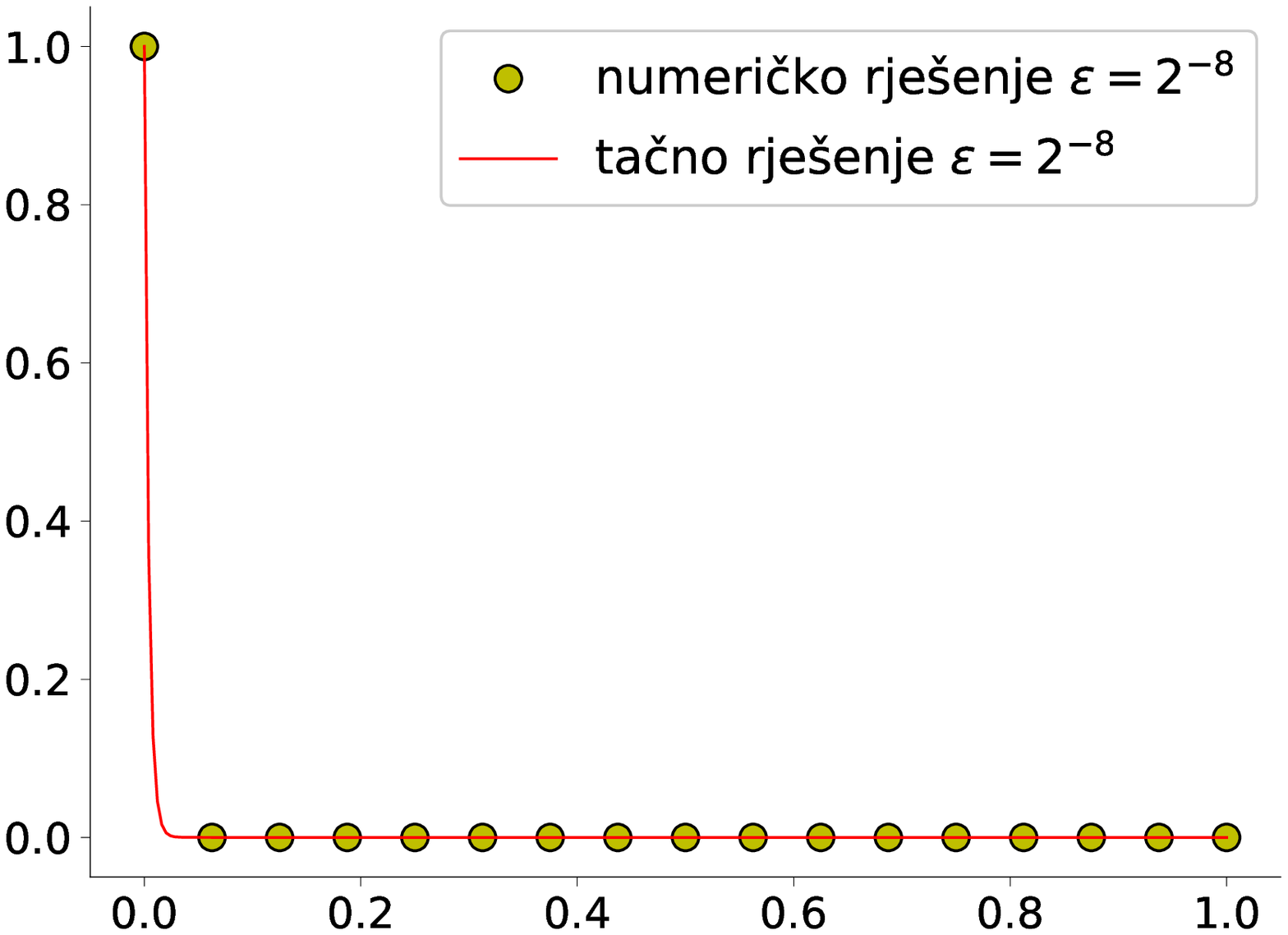}
	\caption{Grafici ta\v cnog i numeri\v ckog rje\v senja na uniformnoj mre\v zi}
	\label{uniformna}
\end{figure}

Radi ilustracije za\v sto uniformne mre\v ze nisu najbolji izbor za rje\v savanje Cauchyjevih problema, \v cija ta\v cna rje\v senja imaju brze promjene, posmatrajmo sliku \ref{uniformna}. Na slici su predstavljeni grafici tri numeri\v cka i tri ta\v cna rje\v senja Cauchyjevog problema \eqref{cauchy1}. 

Gore lijevo su grafici numeri\v ckog i ta\v cnog rje\v senja za vrijednost parametra $\varepsilon=2^{-2}.$ Rubni sloj u okolini ta\v cke   $x=0$ vrlo je slabo izra\v zen, te je lako uo\v citi da su ta\v cke koje predstavljaju numeri\v cko rje\v senje dosta dobro raspore\dj ene na \v citavom grafiku. Dodu\v se ne\v sto rje\dj e u sloju, ali i dalje zadovoljavaju\'ce. 

Na slici gore desno prikazani su grafici obje vrste rje\v senja ali sada za vrijednost parametra $\varepsilon=2^{-4}.$ Ve\'c je bolje izra\v zen sloj i odmah se uo\v cava lo\v sija raspodjela ta\v caka numeri\v ckog rje\v senja. 

Na slici \ref{uniformna} dole su grafici numeri\v ckog i   ta\v cnog rje\v senja ali za vrijednost parametra $\varepsilon=2^{-8}.$ Sloj je sada u veoma uskom dijelu u odnosu na domen, \v sto zna\v ci da su promjene ta\v cnog rje\v senja veoma velike, tj. ta\v cno rje\v senje se veoma brzo mijenja u sloju.  

Treba napomenuti da je kori\v sten isti broj ta\v caka za ra\v cunanje numeri\v ckog rje\v senja u sva tri slu\v caja, te da je rastojanje izme\dj u dvije susjedne ta\v cke mre\v ze ekvidistantno. Sa posljednje slike vidimo da ta\v cke koje predstavljaju numeri\v cko rje\v senje nisu dobro raspore\dj ene, u dijelu grafika koji odgovara sloju nema niti jedna ta\v cke numeri\v ckog rje\v senja, osim ta\v cke koje odgovara po\v cetnom uslovu $y(x_0, \varepsilon)=y_0.$ Ovo nije dobra opcija jer nemamo nikakvu informaciju o rje\v senju u sloju. Da bi se prevazi\v sao ovaj problem mo\v zemo pove\'cati broj ta\v caka, tj. smanjiti rastojanje izme\dj u dvije susjedne ta\v cke, te bi na taj na\v cin podijelili dio domena koji odgovara sloju sa dovoljnim brojem ta\v caka.  Me\dj utim,  ovakav pristup je lo\v s sa prakti\v cne/kompjutacione strane. Naime, sloj je za ovaj problem \v sirine reda $\mathcal{O}(\varepsilon|\ln \varepsilon|).$ Neka je gre\v ska metode reda $\mathcal{O}(h^3),$ \v sto je uobi\v cajeni  red veli\v cine za Runge--Kutta metode drugog reda i neka se zahtijeva da gre\v ska ne bude ve\'ca od $10^{-6},$ ponovimo segment na kojem ra\v cunamo numeri\v cko rje\v senje je $[0,1].$ Sada iz nejednakosti $h^3\leqslant 10^{-6}$ dobijamo da je $h\leqslant 10^{-2},$ drugim rije\v cima segment $[0,1]$ potrebno je podijeliti na 100 podsegmenata odnosno broj ta\v caka je $N=101.$  Ovakavo rezonovanje vrijedi kada ta\v cno rje\v senje nema sloja, ali na\v s Cauchyjev problem ima sloj i tu \v cinjenicu moramo uzeti u obzir. Neka je $\varepsilon=10^{-4}$ i sloj je u ovom slu\v caju \v sirok $\varepsilon|\ln\varepsilon|=10^{-4}|\ln10^{-4}|\approx 9\cdot 10^{-4}.$ Da bi bar jedna ta\v cka mre\v ze bila u sloju mora biti ispunjen uslov $h<0.0009.$ Grubo govore\'ci, broj ta\v caka mre\v ze mora biti ve\'ci od $1111,$ a to je vi\v se od deset puta u odnosu na $N=101$ i to da bi samo jedna ta\v cka  mre\v ze bila u sloju. Naravno, potrebno je mnogo vi\v se ta\v caka od jedne u sloju, pa bi i broj ta\v caka mre\v ze bio mnogo ve\'ci od 1111, odnosno mnogo vi\v se od deset puta  bi trebalo pove\'cati broj ta\v caka  mre\v ze. 

Prethodna diskusija je opisala samo jedan od razloga za\v sto uniformne mre\v ze nisu pogodne za rje\v savanje Cauchyjevih problema koji imaju izra\v zene slojeve. Metode koje se koriste za numeri\v cko rje\v savanje obi\v cnih diferencijalnih jedna\v cina, uobi\v cajeno se poop\v stavaju (kada je to mogu\'ce) i konstrui\v su metode za  numeri\v cko rje\v savanje parcijalnih diferencijalnih jedna\v cina. Ovaj problem sa velikim pove\'canjem broja ta\v cka postao bi i ve\'ci prelaskom na numeri\v cko rje\v savanje parcijalnih diferencijalnih jedna\v cina.

Da se izbjegao prethodno opisani problem potreban je druga\v ciji pristup numeri\v ckom rje\v savanju Cauchyjevih problema  \eqref{problem2}, nego \v sto je pove\'canje broja ta\v caka uniformne mre\v ze. Jedan od mogu\'cih rje\v senje u prevazila\v zenju ovog problema je upotreba mre\v za koje imaju neunifomnu raspodjelu ta\v caka, ovaj pristup pokazao se veoma efikasnim u  numeri\v ckom rje\v savanju rubnih problem.  Vratimo se ponovo na sliku \ref{uniformna} dole.  U slu\v caju kori\v stenja uniformne mre\v ze potrebno je enormno pove\'cati broj ta\v caka da bi ih bilo dovoljno u sloju.  Da bi se izbjeglo pove\'canje broja ta\v caka mre\v ze, a samim tim i nepotrebno pove\'canje vremena ra\v cunanja a i nepotrebno tro\v senje resursa ra\v cunara, potrebno je izvr\v siti druga\v ciju        raspodjelu ta\v caka mre\v  ze. Ovo radimo na sljede\'ci na\v cin: od ukupnog broja ta\v caka mre\v ze $N$ jedan dio ta\v caka koristimo za sloj, dok preostale ta\v cke koristimo za dio mre\v ze koji odgovara domenu van sloja. Rastojanje ta\v caka u sloju je po pravilu manje od rastojanja preostalih ta\v caka koje koristimo van sloja.  

Jedan od na\v cina generisanje ovakve neuniformne (neekvidistantne) mre\v ze vr\v si se generativnom funkcijom $\phi,$ $\phi:[0,1]\mapsto[0,1]$. Ovo je slo\v zena funkcija i sastavljena je najmanje od dvije druge funkcije, jedna slu\v zi za generisanje ta\v caka mre\v ze u sloju, a druga za generisanje ta\v caka van sloja. 

\begin{equation}\label{mreza4}
\phi(\xi, \varepsilon)=
   \begin{cases}
     \phi_1(\xi, \varepsilon),\quad \xi\in[0,\alpha],\\
     \phi_2(\xi, \varepsilon),\quad \xi\in(\alpha,1].
   \end{cases}
\end{equation}

Vrijednost parametra $\alpha$ odre\dj uje koliko \'ce ta\v caka mre\v ze biti u sloju, a koliko van sloja. Kada je poznat analiti\v cki oblik funkcije $\phi$ neuniformnu mre\v zu nije te\v sko generisati. Polazimo od uniformne mreže $\xi_i=ih,\: i=0,1,\ldots,N,\: h=1/N;$  i ta\v cke neuniformne mre\v ze dobijamo na sljede\'ci na\v cin
\begin{equation}\label{mreza5}
     x_i=\phi(\xi_i, \varepsilon),\:i=0,1,\ldots,N.
\end{equation}

Postavlja se pitanje kako konstruisati generativnu funkciju $\phi?$ Polazna ta\v cka za konstruisanje generativne funkcije su procjene izvoda ta\v cnog rje\v senja. Iz Teoreme \ref{teorema2} lako je uo\v citi da u procjeni izvoda ta\v cnog rje\v senja  figuri\v se eksponencijalna funkcija, stoga o\v cekujemo da funkcija kojom generi\v semo ta\v cke mre\v ze u sloju bude neka logaritamska funkcija. Drugi dio generativne funkcije (druga funkcija) je obi\v cno neki polinom kojim generi\v semo ta\v cke mre\v ze van sloja. Ovaj drugi generativne funkcije se bira tako da generativna funkcija bude barem neprekidna na \v citavom domenu.  

Prvu neuniformnu mre\v zu konstruisao je Bakhvalov \cite{bahvalov1969} upravo koriste\'ci logaritamsku funkciju za generisanje ta\v caka u sloju, dok je za generisanje ta\v caka mre\v ze  van sloja koristio linearnu funkciju.  Nakon toga Vulanovi\'c \cite{vulanovic1986doktorat} pojednostavljuje konstruisanje generativne funckije, va\v zno je napomenuti i veliki  doprinos Liseikina \cite{liseikin1989theadaptive, liseikin2001layer, liseikin2018grid} u konstruisanju mnogih generativnih funkcija odnosno neuniformnih mre\v za.

\begin{figure}[!h]
	\centering 
	\includegraphics[scale=.5]{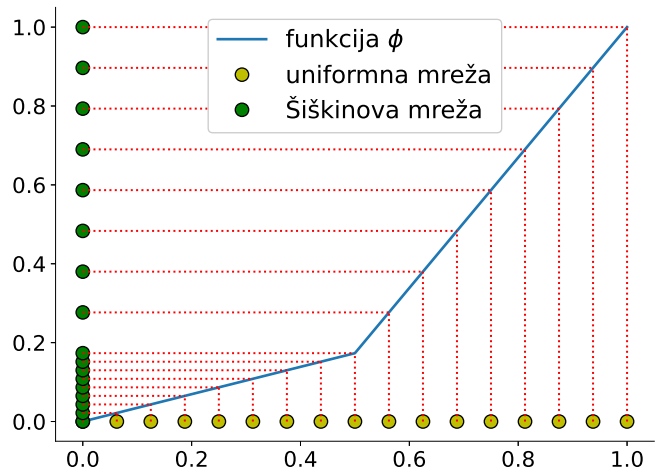}		
	\caption{Generisanje Shishkinove mre\v ze iz uniformne mre\v ze funkcijom $\phi$ }
	\label{mreza3}
\end{figure}

U nastavku ovog rada koristimo najjednostavniju neuniformnu mre\v zu, koju je konstruisao Shishkin \cite{shishkin1988grid}. Ova mre\v za je sastavljena iz dvije uniformne mre\v ze, jedne sa finijom raspodjelom ta\v caka (manje rastojanje izme\dj u ta\v caka) i drugom grubljom (ve\'ce  rastojanje izme\dj u ta\v caka). Shishkin je konstruisao mre\v zu koja slu\v zi za rje\v savanje rubnih problema koji imaju eksponencijalni sloj. Pokazao je da funkcija $\phi_1$ (koja generi\v se ta\v cke u sloju) ne mora biti logaritamskog tipa, nego mnogo jednostavnija linearna fukcija, slika \ref{mreza3}.  Shishkinovu mre\v zu mo\v zemo shvatiti kao dvije uniformne mre\v ze spojene na odgovaraju\'ci na\v cin.  Dobijena mre\v za ima boljih i lo\v sijih osobina u odnosu na mre\v ze koje su konstruisali Bakhvalov i Liseikin.  Jednostavnija je analiza metoda koje koriste Shishkinovu mre\v zu u odnosu na Bakhvalovu mre\v zu i Liseikinove mre\v ze, me\dj utim gre\v ska je ve\'ca kod upotrebe Shishkinove mre\v ze. Va\v zno je napomenuti da je broj problema koji se mogu efikasno rije\v siti upotrebom Shishkinove mre\v ze  mnogo manji od broja problema koji se mogu rije\v siti upotrebom veoma fleksibilnih Liseikinovih mre\v za. Generativna funkcija za Shishkinovu mre\v zu data je sljede\'com formulom

\begin{equation}\label{mreza6}
\phi(\xi,\varepsilon)=
     \begin{cases}
        2\sigma \xi,\quad 0\leqslant \xi \leqslant \alpha,\\
        \sigma+\dfrac{1-\sigma}{1-\alpha}\xi,\quad \alpha <  \xi\leqslant 1, 
     \end{cases}
\end{equation}
gdje je $\alpha\in (0,1),$ $\sigma=\min\{0.5, (n/b)\varepsilon \ln N\},$ $b>0$ i $n\geqslant 1.$  Veli\v cina $\sigma$ je tranziciona ta\v cka mre\v ze ili Shishkinova tranziciona ta\v cka  i ona odre\dj uje mjesto prelaska sa fine na grubu mre\v zu.  Prilikom kori\v stenja Shishkinove mre\v ze za numeri\v cko rje\v savanje rubnih problema, tipi\v cna vrijednost parametra $n$ odgovara stepenu konvergencije metode, dok je parametar $b$ odre\dj en Teoremom o procjeni izvoda i $\alpha=1/2.$

Polaze\'ci sada od generativne funkcije $\phi$ date formulom \eqref{mreza6}, Shishkinovu mre\v zu generi\v semo koriste\'ci formulu \eqref{mreza5}. Na slici \ref{mreza3} dat je grafik funkcije $\phi$ (plava boja). U funkciju  $\phi$ uvr\v stavamo vrijednosti $\xi_i=ih,$ $i=0,1,\ldots,N;$ koje odgovaraju ta\v ckama uniformne mre\v ze. Ove su ta\v cke predstavljene na slici \ref{mreza3} \v zutom bojom i nalaze se na $x$--osi. Odgovaraju\'ce vrijednosti  funkcije $\phi(ih)$ predstavljene su zelenom bojom i predstavljene su na $y$--osi, ovo su ta\v cke Shishkinove mre\v ze, tj. $x_i=\phi(ih)=\phi(\xi_i),\: i=0,1,\ldots, N.$


%

\begin{figure}[!h]
	\centering 
	\includegraphics[scale=.4]{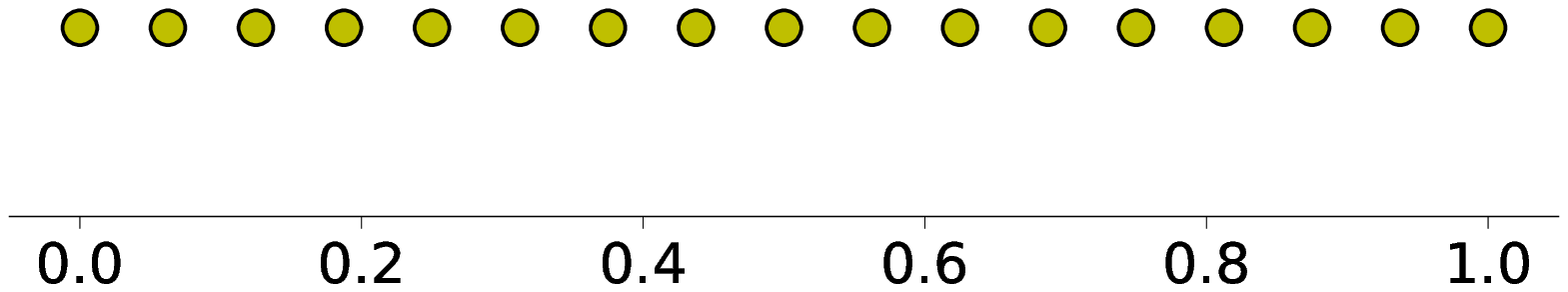}		
	\includegraphics[scale=.4]{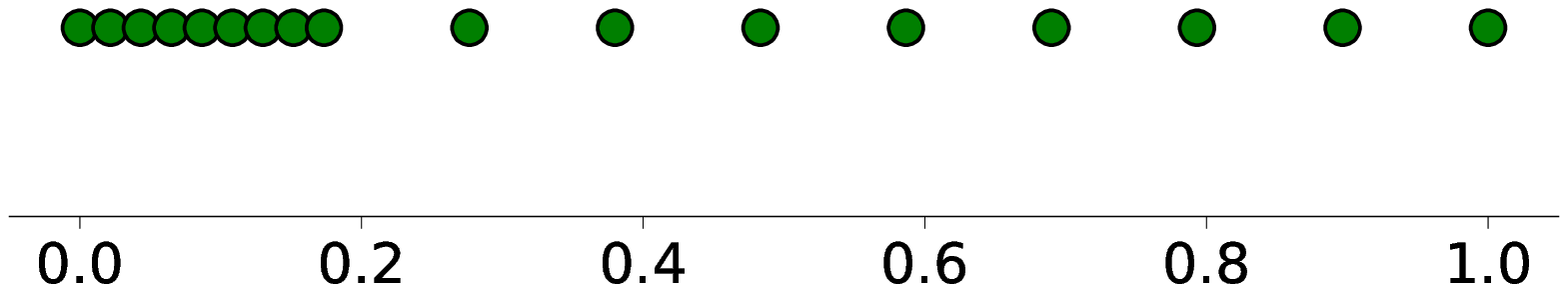}
	\caption{Uniformna (lijevo) i Shishkinova mre\v za (desno)}
	\label{mreza1}
\end{figure}

Na slici \ref{mreza1} predstavljene su ta\v cke uniformne mre\v ze i dobijene ta\v cke Shishkinove mre\v ze. Sa slike \ref{mreza1} (desno) vidi se neuniformna raspodjela ta\v caka mre\v ze, tako\dj e nije te\v sko primjetiti da su ta\v cke mre\v ze kondenzovane u okolini ta\c ke $x=0.$ Ponovimo, u okolini ta\v cke $x=0$ je sloj i da bi se problem koji nastaje zbog brzih  promjena ta\v cnog rje\v senja u sloju rije\v sio na odgovaraju\'ci na\v cin, potrebna je ovakva raspodjela ta\v caka.  
 \\
Vi\v se detalja o slojno--adaptivnim mre\v zama i njihovoj konstrukciji, za razne probleme te razli\v cite vrste slojeva,  mo\v ze se na\'ci u Liseikin i dr.  \cite{liseikin2021numerical}.

\begin{remark}
	Prethodno opisani problem, koji je naveden kao razlog kori\v stenja slojno--adaptivnih mre\v za pri numeri\v ckom rje\v savanje Cauchyjevih problema \eqref{problem2}, u literaturi se naziva i ekonomi\v cnost ra\v cunanja. Osim prevazila\v zenja ovog problema uvo\dj enjem slojno--adaptivnih mre\v za pri numeri\v ckom rje\v savanju problema \eqref{problem2} i sli\v cnih, rje\v sava se i problem stabilnost. Ovo je vrlo ozbiljan problem i prevazilazi okvire ovog rada, stoga \'cemo se stabilnosti samo dota\'ci u numeri\v ckim eksperimentima na kraju rada.. 
\end{remark}

\section{Runge--Kutta metode}
Veoma uspje\v sne u numeri\v ckom rje\v savanju Cauchyjevih problema pokazale su se Runge--Kutta metode \cite{runge, Kutta}. Ovo su jednokora\v cne metode, dakle u Runge--Kutta metodama za ra\v cunanje numeri\v cke vrijednosti u ta\v cki $x_i$ tj. $y_{i+1},$ koristi se samo  jedna predhodno izra\v cunata vrijednost $y_i, $  ona koja odgovara ta\v cki $x_i.$ Ove dvije vrijednosti $x_i$ i $y_i$ koriste se za ra\v cunanje aproksimativne jedne ili vi\v se vrijednosti   funkcije $f.$ Broj ovih aproksimativnih vrijednosti povezan je sa redom ili nivoom Runge--Kutta metoda.

Op\v sti oblik Runge--Kutta metoda $s$--tog reda (ili nivoa) je
\begin{equation}\label{rungekutta1}
y_{i+1}=y_i+h_i\sum_{j=1}^{s}{b_jk_j},
\end{equation}
gdje je 
\begin{equation}\label{rungekutta2}
k_j=f\left(x_i+c_jh_i, y_i+h_i\sum_{q=1}^{s}{a_{j,q}k_q}\right),\: j=1,2,\ldots,s,
\end{equation}
dok su $y_i$ i $y_{i+1}$ aproksimativne (ili numeri\v cke) vrijednosti koje odgovaraju \v cvorovima $x_i$ i $x_{i+1},$ respektivno.

Koeficijente koji se pojavljuju u formulama \eqref{rungekutta1} i \eqref{rungekutta2} predstavljeni su u tabeli \ref{tabela1}, koja predstavlja punu ili implicitnu Runge--Kutta metodu $s$--tog reda. Osim implicitnih postoje i eksplicitne Runge--Kutta metode.  Izme\dj u ove dvije vrste metoda razlika je u sljede\'cem: kod eksplicitnih metoda veli\v cine $k_j,\: j=1,\ldots,s,$ ra\v cunamo koriste\'ci samo prethodno izra\v cunate vrijednosti $k_1,\ldots, k_{j-1},$ za razliku od implicitnih gdje se za ra\v cunanje veli\v cine $k_j$ mogu  koristiti sve veli\v cine $k_1,\ldots, k_s.$  Razlika izme\dj u eksplicitnih i implicitnih metoda dobro je obja\v snjena u literaturi. Uobi\v cajeno, implicitne metode su komplikovanije, samim tim su te\v ze za kori\v stenje. U slu\v caju kada je $f$ nelinearna funkcija potrebno je u svakom koraku rije\v siti nelinearni sistem jedna\v cina. Ali isto tako poznato je da su bolje od eksplicitnih metoda po pitanju stabilnosti.   


\begin{table*}[!h]\centering	
	\begin{tabular}{c| cc cc }
		
		$c_1$& $a_{1,1}$ & $a_{1,2}$ & $\ldots$ & $a_{1,s}$ \\
		$c_2$& $a_{2,1}$ & $a_{2,2}$ & $\ldots$ & $a_{2,s}$ \\
		     $\vdots$&$\vdots$&$\vdots$&&$\vdots$\\
		$c_s$& $a_{s,1}$ & $a_{s,2}$ & $\ldots$ & $a_{s,s}$ \\		
		\bottomrule
	      	& $b_1$&$b_2$& $\ldots$ & $b_s$  
	\end{tabular}
	\caption{Butcherov niz za punu (implicitnu) Runge--Kutta metodu}
	\label{tabela1}
\end{table*}

U tabeli \ref{tabela3} dati su koeficijenti eksplicitne Runge--Kutta metode $s$--tog reda. 

\begin{table*}[!h]\centering	
	\begin{tabular}{c| cc cc c}
		
		$0$& $0$ & $0$ & $\ldots$ &$0$& $0$ \\
		$c_2$& $a_{2,1}$ & $0$ & $\ldots$ & $0$ & $0$\\
		$c_3$& $a_{3,1}$ & $a_{3,2}$ & $\ldots$ & $0$ & $0$\\
		$\vdots$&$\vdots$&$\vdots$&&$\vdots$&\\
		$c_s$& $a_{s,1}$ & $a_{s,2}$ & $\ldots$ & $a_{s,s-1}$& $0$ \\		
		\bottomrule
		& $b_1$&$b_2$& $\ldots$ &$b_{s-1}$& $b_s$  
	\end{tabular}
	\caption{Butcherov niz za eksplicitnu Runge--Kutta metodu}
	\label{tabela3}
\end{table*}

Isto tako  u literaturi pokazano je kako se ra\v cunaju koeficijenti iz tabele \ref{tabela1} ili \ref{tabela3} na vi\v se na\v cina.
\paragraph{Eksplicitne Runge--Kutta metode drugog reda} Sada \'cemo izra\v cunati koeficijente za ekplicitne Runge--Kutta metode drugog reda. Odgovaraju\'ci koeficijenti dati su tabeli \ref{tabela4}.

\begin{table*}[!h]\centering	
	\begin{tabular}{c| cc }
		
		$0$& $0$ & $0$   \\
		$c_2$& $a_{2,1}$ & $0$  \\
		\bottomrule
		& $b_1$&$b_2$  
	\end{tabular}
	\caption{Butcherov niz za eksplicitnu Runge--Kutta metodu reda $s=2$}
	\label{tabela4}
\end{table*}

Numeri\v cku vrijednost $y_{i+1}$ koja odgovara ta\v cki mre\v ze $x_{i+1}$ ta\v cnog rje\v senja $y$ ra\v cunamo po formuli
\begin{equation}
y_{i+1}=y_i+h_i(b_1k_1+b_2k_2),
\end{equation}
gdje su 
\[k_1=f(x_i+c_1h_i, y_i+h_i a_{1,1}k_1)\]
i
\[k_2=f(x_i+c_2h_i, y_i+h_i a_{2,1} k_1).\]
Kako je iz tabele \ref{tabela4}, $c_1=0$ i $a_{1,1}=0,$ to dobijamo
\begin{align}
   y_{i+1}=y_i+h_i\left[ b_1f(x_i,y_i)+ b_2 f(x_i+c_2h_i, y_i+h_i a_{2,1}f(x_i,y_i))  \right].
\end{align}
Koeficijente $b_1,\, b_2,\, c_2$ i $a_{2,1}$ odre\dj ujemo upore\dj uju\'ci Taylorove razvoje za $y_{i+1}$ i $y(x_{i+1}).$ Vrijedi 

\begin{align}\label{razvoj1}
y_{i+1}&= y_i+h_ib_1f(x_i,y_i)+h_ib_2\left[  f(x_i,y_i) + f_x(x_i,y_i)c_2 h_i + f_y(x_i,y_i)h_ia_{2,1}f(x_iy_i)+\mathcal{O}(h^2_i)  \right]
\end{align}
i
\begin{align}\label{razvoj2}
y(x+1)&=y(x_i)+ y'(x_i)h_i+\frac{y''(x_i)}{2}h^2_i+\mathcal{O}(h^3_i)\nonumber\\
      &=y(x_i)+f(x_i,y_(x_i))h_i+\frac{1}{2}\left[ f_x(x_i,y(x_i)) + f_y(x_i,y(x_i))f(x_i,y(x_i))\right]h^2_i+\mathcal{O}(h^3_i).
\end{align}
Upore\dj uju\'ci sada izraze iz \eqref{razvoj1} koji su uz $h_i$ sa izrazima iz \eqref{razvoj2} tako\dj e uz $h_i,$ te ovaj postupak ponovimo i za $h^2_i,$ dobijamo sistem   

\begin{equation}\label{sistem1}
\begin{cases}
      b_1+b_2=1\\
      b_2c_2=\dfrac{1}{2}\\
      b_2a_{2,1}=\dfrac{1}{2}.
\end{cases}
\end{equation}
Sistem \eqref{sistem1} je nelineran sa 3 jedna\v cine i 4 nepoznate i  ne  mo\v  zemo ga jednozna\v cno rije\v siti. Sljede\'ca rje\v senja, odnosno metode se naj\v ce\v s\'ce koriste:

\begin{enumerate}[(a)]
	\item  \ \\
	
	\begin{equation}
	   \begin{cases} \label{RG1}
              b_1=b_2=\dfrac{1}{2},\: c_2=a_{2,1}=1 ,\\
              k_1=f(x_i,y_i),\: k_2=f(x_i+h_i, y_i+k_1),\\
              y_{i+1}=y_i+\dfrac{1}{2}h_i(k_1+k_2),
       \end{cases}
     \end{equation}
     	\item  \ \\
     
     \[\begin{cases} 
     b_1=\dfrac{1}{4},\: b_2=\dfrac{3}{4},\: c_2=a_{2,1}=\dfrac{2}{3}, \\
     k_1=f(x_i,y_i),\: k_2=f(x_i+\frac{2}{3}h_i, y_i+\frac{2}{3}k_1),\\
     y_{i+1}=y_i+\dfrac{1}{4}h_i(k_1+3k_2),
     \end{cases}
     \]
     
     \item  \ \\
     
     \[\begin{cases} 
     b_1=0,\: b_2=1,\: c_2=a_{2,1}=\dfrac{1}{2} ,\\
     k_1=f(x_i,y_i),\: k_2=f(x_i+\frac{1}{2}h_i, y_i+\frac{1}{2}k_1),\\
     y_{i+1}=y_i+ h_i k_2.
         \end{cases}
     \]
\end{enumerate}

\paragraph{Eksplicitne Runge--Kutta metode tre\'ceg reda} Navedeni postupak mo\v zemo iskoristiti i za dobijanje metode Runge--Kuta metoda vi\v seg reda od dva. Za dobijanje metoda tre\'ceg reda, potrebno je samo u formulama \eqref{razvoj1} i \eqref{razvoj2} u razvoj uklju\v citi i druge izvode i iskoristiti izraze uz $h^3_i.$ Na taj na\v cin neke od formula koje dobijamo su

\begin{enumerate}[(a)]
	\item  \ \\
	
	\begin{equation}
	\begin{cases} \label{RG2}
	k_1=f(x_i,y_i),\: k_2=f(x_i+\frac{1}{2}h_i, y_i+\frac{1}{2}k_1),\: k_3=f(x_i+\frac{3}{4}h_i, y_i+\frac{3}{4}k_2 ), \\
	y_{i+1}=y_i+\dfrac{1}{9}h_i( 2k_1+3k_2+4k_3),
	\end{cases}
	\end{equation} 
	\item  \ \\
	
	\[\begin{cases} 
	k_1=f(x_i,y_i),\: k_2=f(x_i+\frac{1}{2}h_i, y_i+\frac{1}{2}k_1),\: k_3=f(x_i+h_i, y_i-k_1+2k_2),\\
	y_{i+1}=y_i+\dfrac{1}{6}h_i(k_1+4k_2+k_3).
	\end{cases}
	\]
\end{enumerate}
Ostale eksplicitne sheme vi\v seg reda dobijaju se na analogan na\v cin.

\paragraph{Implicitna Runge--Kutta metoda drugog reda} Nije te\v sko primjetiti da kod eksplicitnih Runge-Kutta metoda, koeficijente sa ve\'cim indeksom ra\v cunamo koriste\'ci koeficijente sa manjim indeksom, npr. kod eksplicitne Runge--Kutta  metode tre\'ceg reda, $k_2$ ra\v cunamo koriste\'ci $k_1,$ dok $k_3$ ra\v cunamo preko $k_1$ i $k_2.$  Dakle, za ra\v cunanje nekog koeficijenta koristimo već izra\v cunate vrijednosti drugih koeficijenata. Ovo nije slu\v caj sa implicitnim metodama, kao \v sto \'cemo vidjeti na sljede\'cem primjeru. Sljede\'cu implicitnu metodu 

\begin{equation}\label{implicitnaShema1} 
    \begin{cases}
     k_1=f(x_i+(\frac{1}{2}-\gamma)h_i, y_i+\frac{1}{4}h_ik_1+ (\frac{1}{4}-\gamma)h_ik_2) \\
     k_2=f(x_i+(\frac{1}{2}+\gamma)h_i, y_i+(\frac{1}{4}+\gamma)h_ik_1+ \frac{1}{4}h_ik_2) \\       
     y_{i+1}=y_i+\frac{1}{2}h_i(k_1+k_2),
    \end{cases}
\end{equation}
gdje je $\gamma=\frac{\sqrt{3}}{6},$ koristimo u numeri\v ckim eksperimentima. Iz formule \eqref{implicitnaShema1} vidimo da se u izrazu za ra\v cunanje koeficijenta $k_1$ pojavljuje koeficijent $k_2$ i obrnuto. Ovo je razlog za\v sto se ovakve metode nazivaju implicitnim. 

Da bi izra\v cunali numeri\v cku vrijednost $y_{i+1}$ potrebno je da znamo $y_i,\,k_1$ i $k_2,$ zbog toga moramo u svakom koraku rije\v siti sistem 
\begin{equation}\label{implicitnaShema2}
\begin{cases} k_1=f(x_i+(\frac{1}{2}-\gamma)h_i, y_i+\frac{1}{4}h_ik_1+ (\frac{1}{4}-\gamma)h_ik_2) \\ k_2=f(x_i+(\frac{1}{2}+\gamma)h_i, y_i+(\frac{1}{4}+\gamma)h_ik_1+  \frac{1}{4}h_ik_2),
\end{cases}
\end{equation}
po nepoznatim $k_1$ i $k_2.$
U zavisnosti od funkcije $f$ ovo je linearni ili nelinearni sistem od dvije jedna\v cine. U slu\v caju kada je funkcija $f$ linearna, kao \v sto je slu\v caj u \eqref{problem3},   ovaj sistem (koji postaje  sistem dvije linearne algebarske jedna\v cine) mo\v zemo simboli\v cki rije\v siti po $k_1$ i $k_2,$ te ova rje\v senja uvrstiti u izraz $y_{i+1}=y_i+\frac{1}{2}\left(k_1+k_2\right).$ Kada je $f$ nelinearna funkcija situacija je dosta komplikovanija, tada je potrebno u svakom koraku rje\v savati nelinearni sistem.

Ograni\v ci\'cemo se samo na linearni slu\v caj \eqref{problem3}, sada sistem \eqref{implicitnaShema2} poprima novi oblik
\begin{equation}
\begin{cases}
 k_1=p_i^{(1)}\left[ y_i+\frac{1}{4}h_ik_1+\left(\frac{1}{4}-\gamma\right)h_ik_2 \right]+q_i^{(1)} \\
 k_2=p_i^{(2)}\left[ y_i+\left( \frac{1}{4}+\gamma  \right)h_ik_1+\frac{1}{4}h_ik_2  \right]+q_i^{(2)},
\end{cases}
\end{equation}
odnosno

\begin{equation}
   \begin{cases}
   \left(1-\frac{1}{4}p_i^{(1)}h_i\right)k_1-p_i^{(i)}\left(\frac{1}{4}-\gamma\right)h_ik_2=p_i^{(1)}y_i+q_i^{(i)} \\
   -p_i^{(2)}\left(\frac{1}{4}+\gamma\right)h_ik_1+\left(1-\frac{1}{4}p_i^{(2)}h_i\right)k_2=p_i^{(2)}y_i+q_i^{(2)},
   \end{cases}
\end{equation}
gdje su $p_i^{(1)}=p(x_i+\left(\frac{1}{2}-\gamma\right)h_i),$ $p_i^{(2)}=p(x_i+\left( \frac{1}{2}+\gamma \right)h_i),$ analogno vrijedi i za $q_i^{(1)}$ i $q_i^{(2)}.$ Rje\v savaju\'ci prethodni sistem po $k_1$ i $k_2,$ te uvr\v stavaju\'ci dobijene izraze u tre\'ci izraz iz \eqref{implicitnaShema1} dobijamo

\begin{equation}\label{RG3}
y_{i+1}=y_i+\frac{1}{2}h_i\frac{ \left(p_i^{(1)}y_i+q_i^{(1)}\right)\left( 1+p_i^{(2)}\gamma h_i\right) 
	              + \left(p_i^{(2)}y_i+q_i^{(2)}\right)\left( 1-p_i^{(1)}\gamma h_i\right)}
               {\left(1-\frac{1}{4}p_i^{(1)}h_i\right)\left(1-\frac{1}{4}p_i^{(2)}h_i\right)
                -p_i^{(1)}p_2^{(2)}\left(\frac{1}{16}-\gamma^2\right)h_i^2  },\: i=1,2,\ldots,N-1.
\end{equation} 
Prethodna formula predstavlja implicitnu metodu \eqref{implicitnaShema1} u slu\v caju kada je $f$ linearna funkcija ($f(x,y)= p(x)y+q(x)$), tj. kada diferencijalna jedna\v cina linearna.

\section{Numeri\v cki ekperimenti}

U ovoj sekciji bi\'ce testirane metode date u prethodnom dijelu rada na Shishkinovoj mre\v zi koja je generisana funkcijom \eqref{mreza6}.  Veli\v cine koje \'cemo ra\v cunati u ovim numeri\v ckim ekperimentima  su vrijednost gre\v ske $E_N$ i brzina konvergencije Ord.  Ra\v cunamo ih po formulama  

\begin{equation}\label{en}
    E_N=\max_{0\leqslant i\leqslant N}\left| y(x_i)-y^{N}_{i}\right|,
\end{equation}
\begin{equation}\label{ord}
    \ord =\frac{\ln E_N-\ln E_{2N}}{\ln(2k/(k+1))},
\end{equation}
gdje je $N=2^k,$ dok $y(x_i)$ predstavlja vrijednost ta\v cnog rje\v senja Cauchyjevog problema u ta\v cki mre\v ze $x_i$, a $y^{N}_{i}$ je vrijednost odgovaraju\'ceg numeri\v ckog rje\v senja, tako\dj e u ta\v cki mre\v ze $x_i,$ ali koje je izra\v cunato na mre\v zi sa $N+1$ ta\v ckom.

\begin{example}
   Dat je Cauchyjev problem 
   \begin{equation*} 
     \begin{cases}
        \varepsilon y'=-xy+\varepsilon+e^{-x/\varepsilon}+x\left( x-e^{-x/\varepsilon }+1 \right),    \quad 0\leqslant x ,\\
        y(0)=0.
     \end{cases}
   \end{equation*}
   Ta\v cno rje\v senje ovog Cauchyjevog problema je $y(x)=x-e^{-x/\varepsilon}+1.$

   Mre\v zu za ra\v cunanje numeri\v ckog rje\v senja generi\v semo koriste\'ci generativnu funkciju \eqref{mreza3}. Vrijednosti parametara, koji su kori\v steni u numeri\v ckim eksperimentima, su $n=2,\, b=1$ i $\alpha=1/2.$ Vrijednosti perturbacionog parametra $\varepsilon$ i broja ta\v caka $N$ dati su u tabelama.

   U tabelama \eqref{tabelaExp1}, \eqref{tabelaExp2}, \eqref{tabelaExp3}, predstavljene su vrijednosti $E_n$ i $\ord,$ za Runge--Kutta metode  \eqref{RG1}, \eqref{RG2} i \eqref{RG3}, respektivno.

\end{example}

\begin{table*}[!h]\centering	\tiny
	\begin{tabular}{c cc cc cc cc cc }
		\toprule 
		
		\multicolumn{1}{c}{} 
		& \multicolumn{2}{c}{$\varepsilon=2^{-2}$} 
		& \multicolumn{2}{c}{$\varepsilon=2^{-4}$} 
		& \multicolumn{2}{c}{$\varepsilon=2^{-6}$ }  
		& \multicolumn{2}{c}{$\varepsilon=2^{-8}$  }
		& \multicolumn{2}{c}{$\varepsilon=2^{-10}$  } \\
		\cmidrule(r{\tabcolsep}){2-3} 
		\cmidrule(r{\tabcolsep}){4-5} 
		\cmidrule(r{\tabcolsep}){6-7}
		\cmidrule(r{\tabcolsep}){8-9}
		\cmidrule(r{\tabcolsep}){10-11}
		
		\multicolumn{1}{c}{$N$}
		& \multicolumn{1}{c}{$E_N$}
		& \multicolumn{1}{c}{$\ord$}		
		& \multicolumn{1}{c}{$E_N$}
		& \multicolumn{1}{c}{$\ord$}
		& \multicolumn{1}{c}{$E_N$}
		& \multicolumn{1}{c}{$\ord$}
		& \multicolumn{1}{c}{$E_N$}
		& \multicolumn{1}{c}{$\ord$}
		& \multicolumn{1}{c}{$E_N$}
		& \multicolumn{1}{c}{$\ord$}\\
		\midrule
		
		$2^{10}$   &1.49e-06 &2.00 &1.03e-05  &2.00  &1.31e-05  &2.00  &2.78e-05  &2.68    &3.00e-04  &3.30  \\
		$2^{11}$   &4.51e-07 &2.00 &3.13e-06  &2.00  &3.96e-06  &2.00  &5.61e-06  &2.38    &4.09e-05 &3.20 \\
		$2^{12}$   &1.34e-07 &2.00 &9.32e-07  &2.00  &1.81e-06  &2.00  &1.33e-06  &2.04    &5.88e-06 &2.98   \\
		$2^{13}$   &3.94e-08 &2.00 &2.74e-07  &2.00  &3.46e-07  &2.00  &3.81e-07  &2.00    &9.45e-07 &2.69  \\
		$2^{14}$   &1.14e-08 &2.00 &7.93e-08  &2.00  &1.00e-07  &2.00  &1.10e-07  &2.00    &1.78e-07 &2.42  \\
		$2^{15}$   &3.28e-09 &2.00 &2.28e-08  &2.00  &2.88e-08  &2.00  &3.17e-08  &2.00    &3.94e-08 &2.22  \\
		$2^{16}$   &9.34e-10 &2.00 &6.48e-09  &2.00  &8.20e-09  &2.00  &9.01e-09  &2.00    &9.76e-09 &2.07  \\
		$2^{17}$   &2.64e-10 & - &1.83e-09  &-  &2.31e-09  &-  &2.54e-09  &-    &2.64e-09 & -  \\

		\bottomrule
	\end{tabular}
	\caption{Vrijednosti greške $E_N$  i brzine konvergencije $\ord$ za razli\v cite vrijednosti $N$ i $\varepsilon$}
	\label{tabelaExp1}
\end{table*}

\newpage

\begin{table*}[!h]\centering	\tiny
	\begin{tabular}{c cc cc cc cc cc }
		\toprule 
		
		\multicolumn{1}{c}{} 
		& \multicolumn{2}{c}{$\varepsilon=2^{-2}$} 
		& \multicolumn{2}{c}{$\varepsilon=2^{-4}$} 
		& \multicolumn{2}{c}{$\varepsilon=2^{-6}$ }  
		& \multicolumn{2}{c}{$\varepsilon=2^{-8}$  }
		& \multicolumn{2}{c}{$\varepsilon=2^{-10}$  } \\
		\cmidrule(r{\tabcolsep}){2-3} 
		\cmidrule(r{\tabcolsep}){4-5} 
		\cmidrule(r{\tabcolsep}){6-7}
		\cmidrule(r{\tabcolsep}){8-9}
		\cmidrule(r{\tabcolsep}){10-11}
		
		\multicolumn{1}{c}{$N$}
		& \multicolumn{1}{c}{$E_N$}
		& \multicolumn{1}{c}{$\ord$}		
		& \multicolumn{1}{c}{$E_N$}
		& \multicolumn{1}{c}{$\ord$}
		& \multicolumn{1}{c}{$E_N$}
		& \multicolumn{1}{c}{$\ord$}
		& \multicolumn{1}{c}{$E_N$}
		& \multicolumn{1}{c}{$\ord$}
		& \multicolumn{1}{c}{$E_N$}
		& \multicolumn{1}{c}{$\ord$}\\
		\midrule
		

		$2^{10}$   &2.34e-09  &3.00  &5.40e-09   &3.00   &2.23e-09   &3.68   &3.30e-09   &3.04 &2.01e-09   &3.18   \\
		$2^{11}$   &3.86e-10  &3.00  &8.97e-10   &3.00   &2.46e-10   &3.00   &5.34e-10   &3.02 &3.00e-10   &3.10   \\
		$2^{12}$   &6.27e-11  &3.00  &1.45e-10   &3.00   &3.99e-11   &3.00   &8.54e-11   &3.01 &4.56e-11   &3.06   \\
		$2^{13}$   &9.96e-12  &3.00  &2.30e-11   &3.00   &6.33e-12   &3.00   &1.34e-11   &3.00 &6.98e-12   &3.03   \\
		$2^{14}$   &1.55e-12  &3.00  &3.60e-12   &2.99   &9.88e-13   &2.96   &2.09e-12   &2.99 &1.06e-12   &3.01   \\
		$2^{15}$   &2.38e-13  &3.03  &5.55e-13   &2.90   &1.55e-13   &2.61   &3.23e-13   &3.06 &1.62e-13   &2.62   \\
		$2^{16}$   &3.53e-14  &-3.22 &8.94e-14   &2.53   &2.99e-14   &-0.65  &4.70e-14   &-0.21&3.13e-14   &-0.59   \\
		$2^{17}$   &2.70e-13  & -    &1.79e-14   & -     &4.52e-14   & -     &5.39e-14   & -   &4.55e-14   & -  \\

		\bottomrule
	\end{tabular}
	\caption{Vrijednosti greške $E_N$  i brzine konvergencije  $\ord$ za razli\v cite vrijednosti $N$ i $\varepsilon$ }
	\label{tabelaExp2}
\end{table*}

\begin{table*}[!h]\centering	\tiny
	\begin{tabular}{c cc cc cc cc cc }
		\toprule 
		
		\multicolumn{1}{c}{} 
		& \multicolumn{2}{c}{$\varepsilon=2^{-2}$} 
		& \multicolumn{2}{c}{$\varepsilon=2^{-4}$} 
		& \multicolumn{2}{c}{$\varepsilon=2^{-6}$ }  
		& \multicolumn{2}{c}{$\varepsilon=2^{-8}$  }
		& \multicolumn{2}{c}{$\varepsilon=2^{-10}$  } \\
		\cmidrule(r{\tabcolsep}){2-3} 
		\cmidrule(r{\tabcolsep}){4-5} 
		\cmidrule(r{\tabcolsep}){6-7}
		\cmidrule(r{\tabcolsep}){8-9}
		\cmidrule(r{\tabcolsep}){10-11}
		
		\multicolumn{1}{c}{$N$}
		& \multicolumn{1}{c}{$E_N$}
		& \multicolumn{1}{c}{$\ord$}		
		& \multicolumn{1}{c}{$E_N$}
		& \multicolumn{1}{c}{$\ord$}
		& \multicolumn{1}{c}{$E_N$}
		& \multicolumn{1}{c}{$\ord$}
		& \multicolumn{1}{c}{$E_N$}
		& \multicolumn{1}{c}{$\ord$}
		& \multicolumn{1}{c}{$E_N$}
		& \multicolumn{1}{c}{$\ord$}\\
		\midrule
		
		$2^{10}$   &5.84e-5  &5.15  &5.19e-5   &5.17   &4.58e-5   &5.19   &4.01e-5   &5.21 &3.48e-5   &5.23   \\
		$2^{11}$   &2.68e-6  &5.19  &2.35e-6   &5.21   &2.05e-6   &5.23   &1.77e-6   &5.25 &1.52e-6   &5.28   \\
		$2^{12}$   &1.15e-7  &5.16  &9.98e-8   &5.18   &8.59e-8   &5.20   &7.32e-8   &5.22 &6.18e-8   &5.25   \\
		$2^{13}$   &4.84e-9  &5.12  &4.15e-9   &5.14   &3.53e-9   &5.16   &2.97e-13  &5.18 &2.47e-9   &5.20   \\
		$2^{14}$   &2.03e-10 &5.08  &1.72e-10  &5.10   &1.44e-10  &5.12   &1.19e-10  &5.14 &9.84e-11  &5.16   \\
		$2^{15}$   &8.50e-13 &5.08  &7.13e-12  &5.10   &5.91e-12  &5.13   &4.83e-12  &5.12 &3.90e-12  &5.23   \\
		$2^{16}$   &3.48e-13 &3.73  &2.88e-13  &2.31   &2.34e-13  &2.01   &1.93e-13  &2.36 &1.45e-13  &1.19   \\
		$2^{17}$   &4.10e-14 & - &6.68e-14  & -  &6.55e-14  & -  &4.32e-14  &- &6.38e-14  & -  \\

		\bottomrule
	\end{tabular}
	\caption{Vrijednosti greške $E_N$  i brzine konvergencije  $\ord$ za razli\v cite vrijednosti $N$ i $\varepsilon$}
	\label{tabelaExp3}
\end{table*}

\begin{figure}[!h]\centering
	\includegraphics[scale=.4]{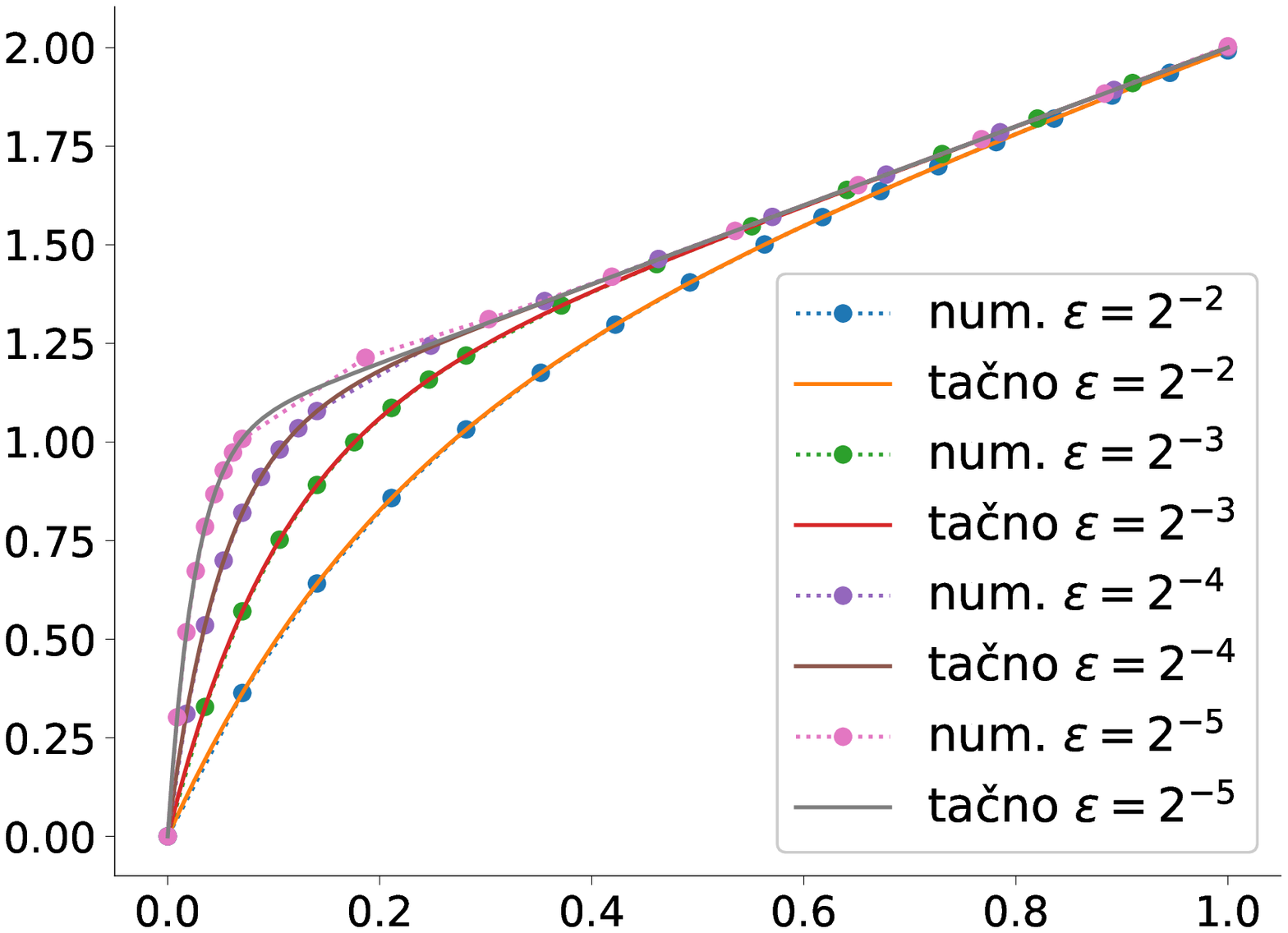}\hspace{-1cm}
	\includegraphics[scale=.4]{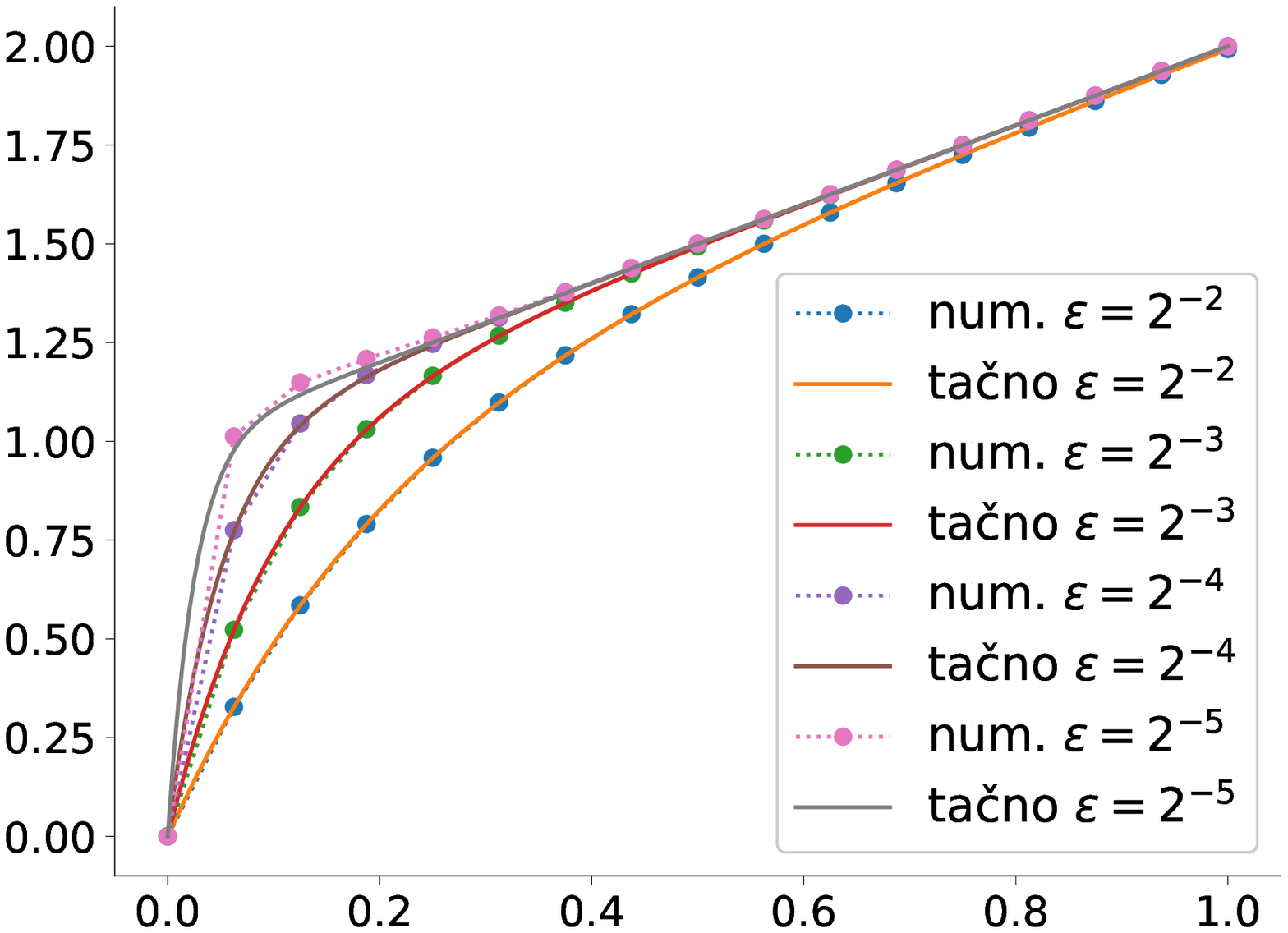}
	
	\includegraphics[scale=.5]{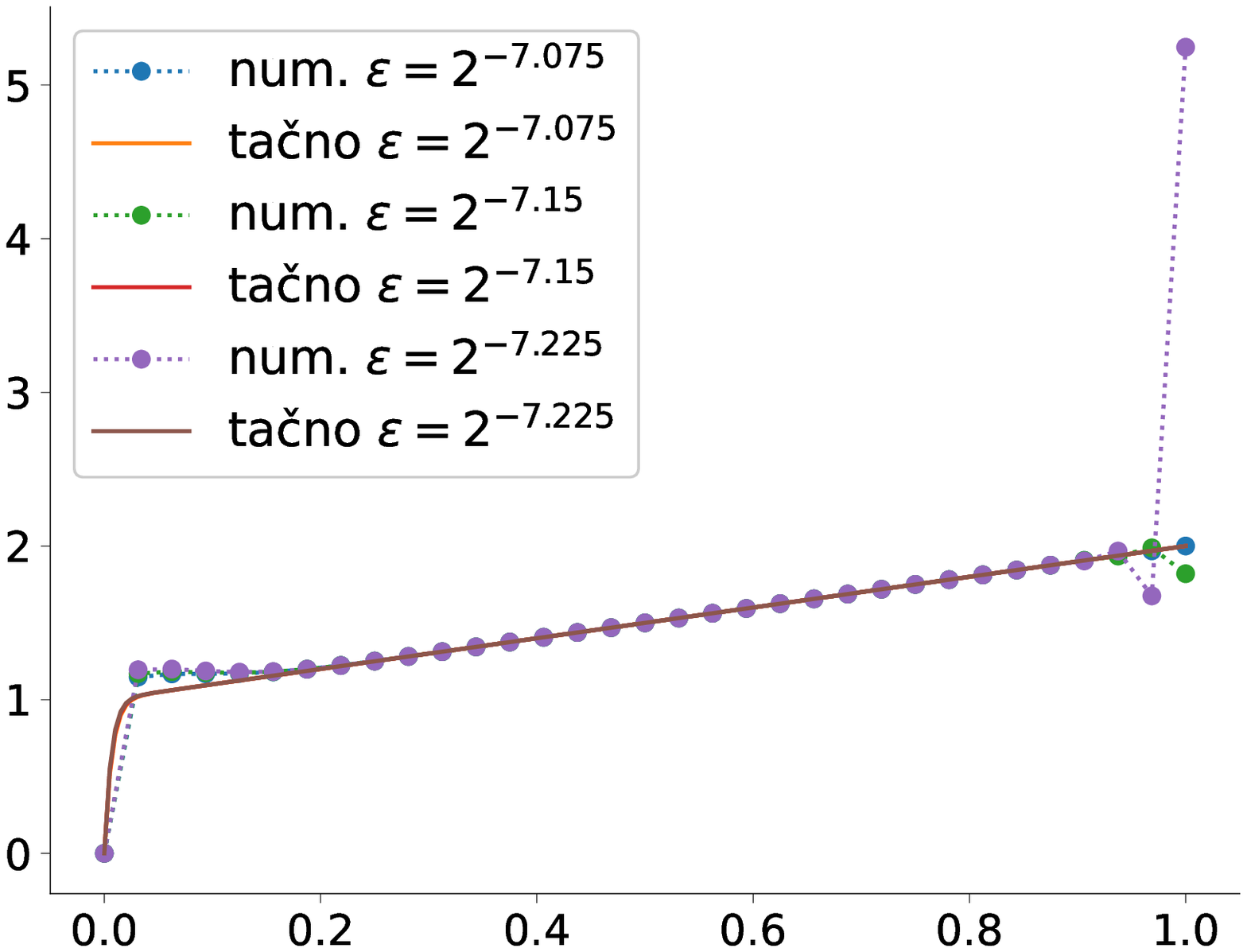}
	\caption{Grafici ta\v cnog i numeri\v ckog rje\v senja na Shishkinovoj (lijevo--gore),  uniformnoj mre\v zi (desno--gore) za $N=32$ i vrijednosti perturbacionog parametra $\varepsilon=2^{-2}, 2^{-3}, 2^{-4}, 2^{-6}$; grafici ta\v cnog i numeri\v ckog rje\v senja na umiformnoj mre\v zi za $N=32$ i $\varepsilon=2^{-7.075},2^{-7.15}, 2^{-7.225}$}
	\label{num1}
\end{figure}

\section{Diskusija i zaklju\v cak}

U ovom radu razmatrano je numeri\v cko rje\v savanje Cauchyjevog problema, \v cije rje\v senje ima izra\v zen sloj. Navedeni su razlozi kori\v stenja slojno--adaptivnih mre\v za pri numeri\v ckom rje\v savanju diferencijalnih jedna\v cina, \v cija ta\v cna rje\v senja imaju brze promjene.  Ova teorija je dobro razvijena za rje\v savanje rubnih problema, i u ovom radu ista metodologija je primjenjena za rje\v savanje Cauchyjevog problema.

Prvo je kori\v stena Runge--Kutta metoda   \eqref{RG1}  na mre\v zi  generisanoj sa \eqref{mreza6}. U tabeli \eqref{tabelaExp1}  su odgovaraju\'ce vrijednosti $E_n$ i $\ord.$ Ova Runge--Kutta metoda je drugog reda, u klasi\v cnoj teoriji (uniformna mre\v za, ta\v cna rje\v senja nemaju izra\v zene brze promjene)  gre\v ska metode je  $\mathcal{O}(h^3),$ odnosno brzina konvergencije je 3. Izra\v cunata vrijednost parametra $\ord$ je 2 ili je bliska ovoj vrijednosti, osim u posljednjoj koloni za vrijednost parametra $\varepsilon=2^{-10},$  gdje vrijednost $\ord$ po\v cinje  sa 3.22 za $N=2^{10}$ i blago se smanjuje do 2.07 za $N=2^{17}.$ Ove vrijednosti nagovje\v stavaju $\varepsilon$--uniformnu konvergenciju, osobinu koja se zahtijeva od metoda koje se koriste za rje\v savanje Cauchyjevih i rubnih problema \v cija rje\v senja imaju izra\v zene brze promjene.  Grubo govore\'ci, metoda ima osobinu $\varepsilon$--uniformne konvergencije, ako vrijednost gre\v ske (dobijena u odgovaraju\'coj normi--uobi\v cajeno maksimum vektorska norma) ne izlazi iz dobijenih teorijskih okvira pri smanjenju perturbacionog parametra $\varepsilon.$ Treba napomenuti da se ra\v cunanje vrijednost parametra $\ord$ blago razlikuje u slu\v caju kori\v stenja uniformne i \v Si\v skinove mre\v ze. 

U tabeli \eqref{tabelaExp2} su vrijednosti $E_n$ i $\ord,$ ali sada je kori\v stena Runge--Kutta metoda \eqref{RG2} na mre\v zi generisanoj sa \eqref{mreza6}. U klasi\v cnoj teoriji vrijednost gre\v ske za ovu metodu je reda $\mathcal{O}(h^4),$ i brzine konvergencije je 4.  Iz prilo\v zene tabele vidimo da je vrijednost parametra $\ord
$ 3 ili je bliska ovoj vrijednosti za broj ta\v caka mre\v ze $N=2^{10}$ do $N=2^{15}.$ Pove\'canjem broja ta\v caka $N=2^{16}$ i dalje, dolazi do smanjenja vrijednosti parametra $\ord,$ te \v cak njegova vrijednost postaje i negativna. Ovakvo pona\v sanje metode nije po\v zeljno i mo\v ze se objasniti akumulacijom vrijednosti gre\v ske, koja je prisutna pri numeri\v ckom   rje\v savanju Cachyjevih problema tipa \eqref{problem2} i sli\v cnih.

U posljednjoj tabeli \eqref{tabelaExp3} su vrijednosti dobijene kori\v stenjem \eqref{RG3} i mre\v ze  generisane sa \eqref{mreza6}. Ovo je Runge--Kutta implicitna metoda drugog reda, vrijednost gre\v ske je $\mathcal{O}(h^3),$ odnosno brzina konvergencije je 3. Dobijene vrijednosti su ve\'ce od 5 za broj ta\v caka mre\v ze $N=2^{10}$ do $N=2^{15},$ pove\'canje  broja ta\v caka vrijednost parametra $\ord$ naglo se smanjuje. Ovako velika odstupanja izra\v cunatih vrijednosti od teorijskih nisu ni u ovom slu\v caju po\v zeljna.

Na slici \eqref{num1} gore lijevo, su  grafici ta\v cnih i numeri\v ckih rje\v senja za razli\v cite vrijednosti perturbacionog parametra $\varepsilon.$ Sva numeri\v cka rje\v senja su dobijena kori\v stenjem $N=33$ ta\v cke. Sa grafika lako je uo\v citi da se smanjivanjem parametra $\varepsilon$ sloj su\v zava, odnosno da je promjena ta\v cnog rje\v senja koji odgovara tom dijelu domena sve br\v za. Evidentno, ta\v cke numeri\v ckog rje\v sanja su dobro raspore\dj ne i u sloju, a to se posti\v ze fleksibilnom konstrukcijom mre\v ze. Sa  druge strane na slici \eqref{num1} gore desno, predstavljeni su grafici ta\v cnih i numeri\v ckih rje\v senja. Ovaj put  za ra\v cunanje numeri\v ckih  rje\v senja kori\v stena je uniformna mre\v za. Nije te\v sko uo\v citi lo\v su osobinu kori\v stenja uniformnih mre\v za za numeri\v cko rje\v savanje problema \eqref{problem2} i njemu sli\v cnih. Naime, smanjivanjem parametra $\varepsilon,$ uz kori\v stenje istog broja ta\v caka mre\v ze, sve je manji broj ta\v caka mre\v ze u sloju, odnosno rastojanje izme\dj u dvije ta\v cke numeri\v ckog rje\v senja postaje neprihvatljivo veliko.     

I na kraju, na  posljednjoj slici \eqref{num1} dole, predstavljeni su grafici ta\v cnih  i numeri\v ckih rje\v senja za razli\v cite vrijednosti parametra $\varepsilon$ i broj ta\v cka $N=33.$ Promjena parametra $\varepsilon$ je veoma mala i na grafiku se ne mogu uo\v citi razlike izme\dj u ta\v cnih rje\v senja. Me\dj utim razlika izme\dj u numeri\v ckih rje\v  senja je velika. Na lijevo strani grafika, vidi se da nema dovoljno ta\v caka mre\v ze u sloju i da su ta\v cke numeri\v ckog rje\v senja previ\v se udaljene jedna od druge, te da je vrijednost gre\v ske velika (rastojanje ta\v caka numeri\v ckog rje\v senja od grafika ta\v cnog rje\v senja). Sa desne strane grafika, situacija je lo\v sija, po\v sto se pojavljuju oscilatorna rje\v senja. Plavim ta\v ckama predstavljeno je numeri\v cko rje\v senje za $\varepsilon=2^{-7.075},$ i ne mo\v se se uo\v citi sa grafika, da ovo rje\v senja odstupa od ta\v cnog  u okolini $x=1.$ Zelenim ta\v ckama predstavljeno je numeri\v cko rje\v senje za $\varepsilon=2^{-7.15},$ sa grafika je lako uo\v citi da se posljednja izra\v cunata vrijednost numeri\v ckog rje\v senja (za $x=1$ ) razlikuje od vrijednosti ta\v cnog rje\v senja--posljednja zelena ta\v cka.  I na kraju, posebno je kriti\v cna situacija sa  numeri\v ckim rje\v senjem za $\varepsilon=2^{-7.225},$ koje je predstavljeno ljubi\v castim ta\v ckama.  Ovo numeri\v cko rje\v senje ima veoma izra\v zeno oscilatorno pona\v sanje. Vrijednost gre\v ske je veoma velika i ovakvo numeri\v cko rje\v senje je potpuno neupotrebljivo, a samim tim njegovo pona\v sanje je potpuno neprihvatljivo. Eliminisanje ovakvih oscilatornih rje\v senja je blisko povezano sa ispitivanjem stabilnosti metode.

Na osnovu testiranih primjera, mo\v zemo zaklju\v citi da je uvo\dj enje slojno--adaptivnih dobra  polazna ta\v cka za numeri\v cko rje\v savanje Cauchyjevih problema tipa \eqref{problem2} i sli\v cnih. Danas je ovo podru\v cje  predmet intezivnog istra\v zivanja i ono se odvija u nekoliko smjerova, aproksimacija izvoda prilago\dj ena ovim problemima, razvoj specifi\v cnih metoda ako i modifikacija slojno--adaptivnih mre\v za.


\newpage
\nocite{*}
\printbibliography

\end{document}